\documentclass[12pt,leqno,twoside]{amsart}
\usepackage{amssymb}
\usepackage{latexsym}
\usepackage{verbatim}
\usepackage{epsfig}

\newtheorem{theorem}{Theorem}[section]
\newtheorem{corollary}[theorem]{Corollary}
\newtheorem{lemma}[theorem]{Lemma}
\newtheorem{proposition}[theorem]{Proposition}
\theoremstyle{definition}
\newtheorem{definition}[theorem]{Definition}
\theoremstyle{remark}

\theoremstyle{remark}

\theoremstyle{remark}

\theoremstyle{remark}
\newtheorem{note}[theorem]{\sc Note}
\theoremstyle{remark}

\theoremstyle{remark}

\renewcommand{\Box}{\square}    
\renewcommand{\Bbb}{\mathbb}
\newcommand{\cal}{\mathcal}

\newcommand{\Cone}{{\rm{Cone}}}

\newcommand{\Hd}{{\rm{Hd}}}  
\renewcommand{\rhd}{{\rm{rhd\hspace{1.5pt}}}}
\newcommand{\rHd}{{\rm{rHd\hspace{1.5pt}}}}
\newcommand{\h}{{\rm{ht}}}

\renewcommand{\i}{{\rm{i}}}

\renewcommand{\int}{{\rm{int}}}

\newcommand{\reg}{{\rm{reg}}}

\newcommand{\Sing}{{\rm{Sing\hspace{1pt}}}}

\newcommand{\id}{{\rm{id}}}

\newcommand{\Iso}{{\rm{Iso}}}
\newcommand{\im}{\mathop{\rm{im}}\nolimits}

\renewcommand{\ker}{\mathop{{\rm{ker}}}\nolimits}

\newcommand{\cl}{{\rm{closure}}}

\newcommand{\ity}{{\infty}}

\newcommand{\fin}{\hspace*{\fill}$\Box$}
\newcommand{\m}{\setminus}

\newcommand{\var}{{\rm{var}}}

\newcommand{\cD}{{\cal D}}

\newcommand{\cS}{{\cal S}}

\newcommand{\cW}{{\cal W}}

\newcommand{\cN}{{\cal N}}

\newcommand{\bR}{{\Bbb R}}
\newcommand{\bC}{{\Bbb C}}

\newcommand{\bP}{{\Bbb P}}

\newcommand{\bZ}{{\Bbb Z}}
\newcommand{\bX}{{\Bbb X}}
\newcommand{\bY}{{\Bbb Y}}

\begin{document}

\title[Vanishing cycles of pencils of hypersurfaces]{Vanishing cycles of  pencils of hypersurfaces}
\author{Mihai Tib\u ar}

\address{Math\' ematiques, UMR 8524 CNRS,
Universit\'e des Sciences et Technologies de Lille, \  59655 Villeneuve d'Ascq, France.}
\email{tibar@agat.univ-lille1.fr}
\dedicatory{To Dirk Siersma, on the occasion of his 60th anniversary}
 
\maketitle

\begin{abstract} 
We prove an extended Lefschetz principle for a large class of pencils of hypersurfaces having isolated singularities, possibly in the axis, and show that the module of vanishing cycles is generated by the images of certain variation maps. 
\end{abstract}

\setcounter{section}{0}
\section{Introduction}

In this paper we extend the Lefschetz principle of slicing by pencils to nongeneric pencils of hypersurfaces on singular non-compact spaces. 
  We started to develop this point of view in \cite{Ti-lef} for proving connectivity theorems of Lefschetz type for nongeneric pencils. Here we go further and introduce global variation maps in order to
  control vanishing cycles. As a result, we prove 
 a far reaching extension of the Second Lefschetz Hyperplane Theorem.

To get an idea of the main result, let us first briefly recall the classical Lefschetz Hyperplane Theorems (see also Note \ref{n:hist} for some references).
 For a projective manifold $Y$ and some hyperplane section $Y_0$, the First Lefschetz Theorem tells us that the map (induced by inclusion):
\begin{equation} \label{eq:homol}
 H_j(Y_0, \bZ) \to H_j(Y, \bZ)
\end{equation} 
 is bijective for $j< n-1$ and surjective for $j=n-1$. The kernel of the surjection
 in dimension $n-1$ is described by the second Lefschetz theorem, whenever $Y_0$ is a generic member of a {\em generic pencil}, i.e.
 the pencil has only complex Morse critical points. Loosely speaking, each such critical point produces a local vanishing cycle and those vanishing cycles together generate $\ker (H_j(Y_0, \bZ) \to H_j(Y, \bZ))$.

We consider here a very general situation: a complex analytic space  $X = Y\setminus V$ with arbitrary singularities, where $Y$ is some compact complex space and $V$ is a complex analytic subspace. We also consider more general divisors than hyperplanes, namely {\em pencils of hypersurfaces on $Y$}.
Taking hypersurfaces instead of hyperplanes considerably enlarges the class under study. Let us remark that any singular holomorphic function germ is a local pencil of hypersurfaces (but not of hyperplanes): our approach enfolds the theory of hypersurface singularities started by Milnor \cite{Mi}, in some of its aspects.  In the global setting, the polynomial functions $\bC^n \to \bC$ with isolated singularities represent a distinguished class of pencils of hypersurfaces, thus non-hyperplane pencils, which we shall discuss in \S \ref{s:poly}. For more details on this viewpoint and for examples we refer to \cite{Ti-mero} and  \cite{LT}.

  On the other hand, we weaken the classical genericity condition ``the axis of the pencil is in general position in $Y$'' by allowing that the genericity of the axis fails at a finite number of points.  This conceptual extension is introduced and explained in detail in our paper \cite[\S 2]{Ti-lef}. Pencils
which allow such ``isolated singularities in the axis'' (see \S \ref{ss:pencils} for the definitions)
 are natural to consider since isolated singularities of functions on singular spaces are central objects of study in 
  modern singularity theory. (One may refer to the pioneering work, on stratified Morse theory by Goresky and MacPherson \cite{GM}, and on the topology of functions on singular spaces by L\^e D\~ ung Tr\'ang \cite{Le}.)
  
We show in \S \ref{pencils} that one can define the following {\em variation map} around each critical value $a$ of the pencil:
\[ \var_a : H_*(X_c, (X_a)_\reg) \to H_*(X_c),\]

 and that the module of {\em vanishing cycles} at $X_a$, i.e. the kernel of the surjection similar to (\ref{eq:homol}), is generated by the images of these variation maps. Our variation maps can be viewed as global versions of the local variation maps that one defines in singularity theory, see e.g. \cite{La}, \cite{Lo}, \cite[4.4]{Ti-compo}, \cite[\S 2]{NN2}.
 It is well-known that in case of non-isolated singularities,  the local  variation maps do not exist. This is the main reason why the use of variation maps in our results would not extend to this context. Let us remark that in case of one-dimensional singularities, Siersma \cite{Si} defined other types of variation maps, but their behaviour appears to be much more delicate and has not been exploited yet in the literature.   
   
As for our approach, it starts in the spirit of the Lefschetz method \cite{Lef}, as presented by Thom in his Princeton talk in 1957 and by Andreotti and Frankel in their paper \cite{AF2}. This vein has been exploited in relatively few papers ever since; we may mention the interesting ones by Lamotke \cite{La}, Ch\'eniot \cite{Ch1, Ch} and Eyral \cite{Ey}.  The use, in the statement of our Theorem \ref{t:main}, of the comparison between the general element of the pencil and the axis comes from Lamotke \cite{La} and may evoke Ch\'eniot's statements in {\em loc.cit.}  Our setting being far more general, we follow a different strategy and use in a crucial way specific geometric constructions and results of stratified singularity theory.

 A highly nongeneric situation is encountered when the axis of the pencil is contained in $V$. We show that if $V$ contains a member of the pencil then, surprisingly, Theorem \ref{t:main} and its proof still work, with even less restrictive assumptions. 
  Actually, one of the reasons to study such nongeneric pencils is that the polynomial functions on $\bC^n$ constitute a class of examples.  
We show in \S \ref{s:poly} how Theorem \ref{t:2} can be extended to a polynomial function with isolated singularities, but without any condition on the singularities along the axis (which are the so-called ``singularities at infinity''). We arrive in this way to results on vanishing cycles of polynomials which have been discovered in somewhat different form by Neumann and Norbury \cite{NN1, NN2}, such as an invariant cycle theorem (Corollary \ref{c:poly}).

In \S \ref{s:comments} we compare the assumptions of our main Theorem \ref{t:main} 
to conditions involving the {\em rectified homological depth} (defined by Grothendieck and thoroughly studied in \cite{HL}), showing that the latter are more restrictive. We also point out how, by relaxing the generality of the setting, one may recover several results in the literature.

 This paper is based on our preprint \cite{Ti-lefht} and is a natural continuation of \cite{Ti-lef}. We wish to thank the Newton Institute at Cambridge and the Institute for Advanced Study at Princeton for support
 durning the elaboration of this work. We are also thankful to the anonimous
 referee for his valuable remarks.

\section{Nongeneric pencils and variation maps}\label{pencils}    

 Let $Y$ be a compact complex analytic space and let 
  $V\subset Y$ be a complex analytic subspace such that 
   $X := Y\m V$ is of dimension $n$, $n\ge 2$.  
\subsection{Pencils with singularities in the axis}\label{ss:pencils}
    
 Let us recall some definitions that we already used in \cite{Ti-lef}.
  By {\em pencil} (or meromorphic function) we mean the ratio of two
sections $f$ and $g$ of a holomorphic line bundle $L\to Y$.
This defines a holomorphic function $h:= f/g$ over the complement $Y\m A$ of the axis of the pencil $A:= \{ f=g = 0\}$. A pencil is called {\em generic} with respect to $X$ when $Y$ is embedded in some manifold $Z$ and the pencil extends to one over $Z$ which satisfies the following conditions: the axis $\hat A$ of the extended pencil is nonsingular and
transversal to some Whitney stratification of the pair
$(Y,V)$
and the holomorphic map $h= f/g : Y\m A \to \bP^1$ has only stratified
double points as singularities. Notice that part of those singularities might be on $V$, hence outside $X$. 

 Here we focus on a class of (nongeneric) pencils, namely pencils having at most 
isolated singularities, possibly in the axis. 
Let us first explain what we mean by singularities of a pencil.
 

We define a new space by blowing-up along the base locus $A$.
 The idea of this construction is due to Thom and was used by Andreotti and Frankel \cite{AF2} in case of generic pencils on projective manifolds. So, let:
\[  \bY :=\cl \{ (y,[s:t]) \in Y\times \bP^1 \mid sf(y) - tg(y) =0\}.\]
This is a hypersurface in $Y\times \bP^1$ obtained as a Nash blowing-up of $Y$ along $A$. 
 It is clear that the intersection $\bY \cap (Y\m A)\times \bP^1$ is just the graph of $h$, hence it is isomorphic to $Y\m A$. 
 It also follows that the subset $A\times \bP^1$ is included into $\bY$.

Let us denote $\bX:= \bY \cap (X\times \bP^1)$. Consider the projection $p : \bY \to \bP^1$,
its restriction $p_{|\bX} : \bX \to \bP^1$ and
the projection to the first factor $\sigma : \bY \to Y$.
 Notice that the restriction of $p$ to $\bY\m(A\times \bP^1)$ can be identified with $h$.

Now fix a stratification $\cW$ on $Y$ such that $V$ is a union of strata. The restriction of $\cW$ to the open set $Y\m A$ induces a Whitney 
stratification on $\bY\m(A\times \bP^1)$, via the above mentioned identification.
We then denote by $\cS$ the coarsest Whitney
stratification on $\bY$ which coincides over $\bY\m (A\times\bP^1)$ with
the one induced by $\cW$ on $Y\m A$. This stratification exists within a neighbourhood of $A\times \bP^1$, by usual 
arguments (see e.g. \cite{GLPW}), hence such stratification is well defined on $\bY$.  We call it the
{\em canonical stratification} of $\bY$ generated by the stratification $\cW$ 
of
$Y$. The canonical stratification of $\bX$ will be the restriction
of $\cS$ to $\bX$.

\begin{definition} \label{d:sing} 
 We call the {\em singular locus} of $p$ with respect to $\cS$ the following closed analytic subset of $\bY$:
 \[  \Sing_\cS p := \bigcup_{\cS_\beta \in \cS} \Sing
p_{|\cS_\beta}.\]
We denote by $\Lambda := 
p(\Sing_\cS p)$ the set of {\em critical values} of $p$ with respect to $\cS$.
\end{definition}

Since $p$ is proper 
and since $\cS$ has finitely many strata, it follows that the set $\Lambda$ is a finite set.  
By Thom's Isotopy Lemma \cite{Th}, we get  
 that the maps  $p: \bY\m p^{-1}(\Lambda) \to 
\bP^1\m \Lambda$ and $p_{|\bX}:
\bX \m p^{-1}(\Lambda) \to \bP^1\m \Lambda$ are  stratified
locally trivial fibrations. In particular, $h: Y\m (A\cup h^{-1}(\Lambda)) \to \bP^1\m \Lambda$ is a locally trivial fibration. 
   

\begin{definition}\label{d:nongeneric}
We say that the pencil defined by the meromorphic function $h = f/g$ is a {\em pencil 
with isolated singularities} if $\dim \Sing_\cS p \le 0$.

We shall say that $X$ has the structure of a {\em Lefschetz fibration with isolated singularities} if there exists a pencil on $X$ with isolated singularities.
\end{definition}

 We have pointed out in \cite[\S 2]{Ti-lef} that in case $Y$ is projective, the condition $\dim \Sing_\cS p \le 0$ is equivalent to the following condition: the singularites of the 
function $p$ at the blown-up axis $A\times \bP^1$ are at most isolated.
We have moreover:

\begin{proposition} \cite[Proposition 2.4]{Ti-lef}\label{l:lefsing}
 Let $Y\subset \bP^N$ be a projective variety endowed with some Whitney stratification $\cW$ and let $\hat h=\hat f/\hat g$ define a pencil of hypersurfaces in $\bP^N$ with axis $\hat A$.  Let $S$ denote the set of points on $\hat A \cap Y$ where some member of the pencil is singular or where $\hat A$ is not transversal to $\cW$. If $\dim S\le 0$ and the singular points of
$h: Y\m A\to \bP^1$ with respect to $\cW$ are at most isolated, then $\dim \Sing_\cS p \le 0$.
\fin
\end{proposition}
 
\subsection{Variation maps} \label{var}

 
 We assume that our pencil defined by $h\colon  Y\dashrightarrow \bP^1$ has isolated singularities, as defined in \ref{d:nongeneric}.  
\label{piece:var}
Let us fix some notation. For any $M\subset \bP^1$, we denote $\bY_M := p^{-1}(M)$, $\bX_M := \bX \cap \bY_M$, $Y_M := \sigma(p^{-1}(M))$ and $X_M := X \cap Y_M$.
Let $\Lambda = \{ a_1, \ldots , a_p\}$. We denote by $a_{ij}\in \bY$ some point of $\Sing_\cS p \cap p^{-1}(a_i)$. We then have $\Sing_\cS p = \cup_{i,j}\{a_{ij}\}$. For $c\in \bP^1\m \Lambda$ we say that $\bY_c$, resp. $\bX_c$, is a {\em general fiber} of $p : \bY\to \bP^1$, resp. of $p_{|\bX} :\bX \to \bP^1$.
We say that $Y_c$, resp. $X_c$, is a general member of the pencil on $Y$, resp. on $X$.

At some singularity $a_{ij}$, in local coordinates, we take a ball $B_{ij}$ centered at $a_{ij}$. For a small enough radius of $B_{ij}$,
this is a ``Milnor ball" of the holomorphic function $p$ at $a_{ij}$.  Next we may take a small enough disk $D_i\subset \bP^1$ at $a_i \in \bP^1$, so that $(B_{ij}, D_i)$ is Milnor data for $p$ at $a_{ij}$. Moreover, we may do this for all (finitely many) singularities in the fiber $\bY_{a_i}$, keeping the same disk $D_i$, provided it is small enough.

Now the restriction of $p$ to $\bY_{D_i} \m \cup_j B_{ij}$ is a trivial fibration over $D_i$. One may construct a stratified vector field which trivializes this fibration and such that this vector field is tangent to the boundaries of the balls $\bY_{D_i} \cap \partial \bar B_{ij}$. Using this, we may also construct a geometric monodromy of the fibration $p_| : \bY_{\partial \bar D_i} \to \partial \bar D_i$ over the circle $\bar D_i$, such that this monodromy is the identity on the complement of the balls, $\bY_{\partial \bar D_i} \m \cup_j B_{ij}$. The same is then true, when replacing $\bY_{\partial \bar D_i}$ by $\bX_{\partial \bar D_i}$.

Take some point $c_i \in \partial \bar D_i$. We have the geometric monodromy representation:
\[ \rho_i : \pi_1 (\partial \bar D_i, c_i) \to \Iso (X_{c_i}, X_{c_i}\m \cup_j B_{ij}), \]
where $\Iso (.,.)$ denotes the group of relative isotopy classes of stratified homeomorphisms (which are C$^\ity$ along each stratum).
It follows that the geometric monodromy restricted to $X_{c_i}\m \cup_j B_{ij}$ is the identity.

As shown above, we may identify, in the trivial fibration over $D_i$, the fiber $X_{c_i}\m \cup_j B_{ij}$ to the fiber $X_{a_i}\m \cup_j B_{ij}$. Furthermore, in local coordinates at $a_{ij}$, $X_{a_i}$ is a germ of a complex analytic space; hence, for a small enough ball $B_{ij}$, the set $B_{ij}\cap X_{a_i}\m \cup_j a_{ij}$ retracts to $\partial \bar B_{ij}\cap X_{a_i}$, by the local conical structure of analytic sets \cite{BV}.
  Therefore $X^*_{a_i} := X_{a_i}\m \cup_j a_{ij}$ is homotopy equivalent, by retraction, to $X_{a_i}\m \cup_j B_{ij}$.
   
\smallskip  
\noindent {\bf Notation.}  From now on, we shall freely use  $X^*_{a_i}$ as notation for $X_{c_i}\m \cup_j B_{ij}$ whenever we consider the pair $(X_{c_i}, X^*_{a_i})$, having in mind the homotopy equivalence between the two spaces.

It then follows that the geometric monodromy induces an algebraic monodromy, in any dimension $q$:
\[ \nu_i \colon H_q (X_{c_i}, X^*_{a_i}; \bZ) \to  H_q (X_{c_i}, X^*_{a_i}; \bZ),\]
such that the restriction $\nu_i \colon H_q (X^*_{a_i})\to  H_q(X^*_{a_i})$ is the identity.
 
 Consequently, any relative cycle $\delta \in H_q (X_{c_i}, X^*_{a_i}; \bZ)$ is sent by the morphism $\nu_i - \id$ to an absolute cycle. In this way we define a {\em variation map}, for any $q\ge 0$:  
 
\begin{equation}\label{eq:var}
 \var_i : H_q(X_{c_i}, X^*_{a_i}; \bZ) \to H_q(X_{c_i}; \bZ).
\end{equation} 

This enters, as a diagonal morphism, in the following diagram:
\[ \begin{array}{ccc}
H_q (X_{c_i}) & \stackrel{\nu_i - \id}{\longrightarrow} & H_q (X_{c_i}) \\
\mbox{\tiny $j_*$} \downarrow  & \mbox{\tiny $\var_i$} \nearrow  & \downarrow \mbox{\tiny $j_*$} \\
 H_q(X_{c_i}, X^*_{a_i}) & \stackrel{\nu_i - \id}{\longrightarrow} & H_q(X_{c_i}, X^*_{a_i})
\end{array} \]
where $j_*$ is induced by inclusion.

Variation morphisms enter traditionally in the description of global and local fibrations of holomorphic functions at singular fibers, see e.g. \cite{Mi}, \cite{La}, \cite{Si}, \cite[\S 2]{NN2}.  In dimension 2, already Zariski used $\nu_i - \id$ in his theorem for the fundamental group.
 Ch\'eniot \cite{Ch} also works with a kind of a variation map, different from ours. Our definition is a direct extension of the local variation maps (see e.g. \cite{La, Lo}) to the global setting. 

\section{The Main Theorem}\label{main}

Let us recall the definition of the homological depth of a topological space at a point.
\begin{definition}\label{d:hd}
For a discrete subset $\Phi \subset \bX$, we denote by $\Hd_\Phi \bX$ the {\em homological depth} of $\bX$ at $\Phi$. We say that $\Hd_\Phi \bX \ge q+1$ if, at any point $\alpha\in \Phi$, there is an arbitrarily small neighbourhood
$\cN$ of $\alpha$ such that  $H_i(\cN, \cN\m \{\alpha\})= 0$, for $i\le q$.
\end{definition}

For a manifold $M$, at some point $\alpha$, we have $\Hd_\alpha M \ge \dim_\bR M$. Complex $V$-manifolds are rational homology manifolds.
So the homological depth measures the defect of being a homology manifold (for certain coefficients).
 For stratified complex spaces, Grothendieck \cite{G} introduced the rectified homotopical depth, respectively the {\em rectified homological depth}, denoted  $\rHd$. This were later investigated by Hamm and L\^e \cite{HL}, who proved several of Grothendieck's conjectures regarding them. 
  See Proposition \ref{n:cond} for more details and results involving $\rHd$.

 We may now state our principal result, using the notations in \S \ref{pencils}. The homology is with coefficients in $\bZ$.
  

\begin{theorem}\label{t:main}
  Let $h\colon Y \dashrightarrow \bP^1$ define a Lefschetz fibration on $X = Y\m V$ with isolated singularities (cf Definition \ref{d:nongeneric}). Let the axis A be not included in $V$.
    For some $k\ge 0$, suppose that the following conditions are fulfilled:
    
\smallskip\smallskip
\noindent
{\bf (C1)}  \ \ \ $H_q(X_c, X_c\cap A)=0$  for $q\le k$.\\
{\bf (C2)}  \ \ \ $H_q(X_c, X^*_{a_i})=0$ for $q\le k$ and for all $i$. \\ 
{\bf (C3)}  \ \ \ $\Hd_{\bX \cap \Sing_\cS p} \bX \ge k+3$. 
\smallskip\smallskip

 Then $H_q(X, X_c)=0$ for $q\le k+1$ and the kernel of the surjection $H_{k+1}(X_c) \twoheadrightarrow H_{k+1} (X)$ is generated by the images of the variation maps $\var_i$, for $i=\overline{1,p}$.
\end{theorem}


\begin{note}\label{n:main}  
For the annulation of the relative homology we need in fact a weaker condition than (C3), namely the following:

\smallskip
\noindent
{\bf (C3i)}  \ \ \   $\Hd_{\bX \cap \Sing_\cS p} \bX \ge k+2$. 
\smallskip

This will be clear from the proof, since (C3) is used (with $k+3$) only in Corollary \ref{c:2} and Proposition \ref{p:3}(b). See also Proposition \ref{o:1} for what become conditions (C2) and (C3) in special cases, and  Proposition \ref{n:cond} for comparison to the rectified homological depth condition.
For instance, it is well-known from \cite{HL} that, in case $X$ is a {\em complete intersection}, then $\rHd X\ge \dim_\bC X$. This implies (see e.g. the proof of Proposition \ref{n:cond}) that condition {\bf (C3)} is satisfied in this case for $k\le \dim_\bC X -3$. 
\end{note}


In \S 4, we derive the form of this result in special cases, such as in case $\Sing_\cS p \cap (A\times \bP^1) \cap \bX = \emptyset$ (i.e. ``no singularities in the axis''), in case the Lefschetz structure of the space $X$ is hereditary on slices and also in the complementary case $A\subset V$.
 
\begin{note}\label{n:hist}
During the time, Lefschetz hyperplane theorems have been generalized in several directions, giving rise to an extended literature, which the limited space does not allow us to cite here. May we just refer to Fulton's general overview \cite{Fu}, Lamotke's ``classical'' modern presentation of Lefschetz theorems \cite{La} and to Goresky-MacPherson's book \cite{GM} which covers a lot of material.
    
     On the other hand, the description of the kernel of the surjection stated above has been considered in a few papers only. The most recent results are for generic pencils of hyperplanes on quasi-projective manifolds, by Ch\'eniot \cite{Ch}, and on complements in $\bP^n$ of hypersurfaces with isolated singularities and for higher homotopy groups, by Libgober \cite{Li}. The extension of Theorem \ref{t:main} to homotopy groups is investigated in
the preprints \cite{Ti-lefht, Ti-pre}; see also \cite{CL}.

In \S 4 we compare the conditions of our Theorem \ref{t:main} (and show that they are significantly less restrictive) to the conditions used by some other authors in more particular settings than ours: the rectified homology depth condition used by 
Hamm and L\^e \cite{HL-zar, HL-gen, HL} (see \S \ref{ss:rHd}), respectively conditions used by Ch\' eniot and Eyral \cite{Ch1,Ch, Ey} (see \S \ref{ss:relax}).   
\end{note}
\subsection{Proof of Theorem \ref{t:main}}

 Let $K\subset \bP^1$
be a closed disk with $K\cap \Lambda = \emptyset$ and let $\cD$ denote the closure of its complement in $\bP^1$. We denote by $S:= K\cap \cD$ the common boundary, which is a circle, and take a point $c\in S$.
Then take standard paths $\gamma_i \subset \cD \m \cup_iD_i$ (non self-intersecting, non mutually intersecting) from $c$ to $c_i\in \partial \bar D_i$. The configuration $\cup_i(\bar D_i \cup \gamma_i)$ is a deformation retract of $\cD$.
We shall also identify all fibers $X_{c_i}$ to the fiber $X_c$,  by parallel transport along the paths $\gamma_i$.
   
 We denote  $A':= A\cap X_c$. Since $A\not\subset V$, we have that $A' \not= \emptyset$.

\begin{proposition}\label{p:1}
 If $H_q(X_c, A')=0$ for $q\le k$, then the morphism induced by inclusion:
 \[  H_q(X_\cD, X_c) \stackrel{\mbox{\tiny $\iota_*$}}{\longrightarrow} H_q(X, X_c) \]
 is an isomorphism for $q\le k+1$ and an epimorphism for $q= k+2$. 
\end{proposition}
\begin{proof}
 We claim that, if $H_q(X_c, A')=0$, for $q\le k$, then $H_q(X_S, X_c)=0$  for $q\le k+1$. 
 Note first that $X_S$ is homotopy equivalent to the subset
 $\bX_S \cup ( A'\times K)$ of $\bX_K$.
 Let $I$ and $J$ be two arcs which cover
 $S$. We have the homotopy equivalence $(X_S, 
X_c)\stackrel{\h}{\simeq}
 (\bX_I\cup  (A'\times K) \bigcup \bX_J \cup  (A'\times K), \bX_J \cup
 ( A'\times K))$.
 Then, by excision, we have the isomorphism:
 \[ H_*(X_S, X_c) \simeq H_*(\bX_I\cup  (A'\times K), \bX_{\partial I}\cup  (A'\times K)).\]
Furthermore, we have the homotopy equivalences of pairs:
  $(\bX_I\cup  (A'\times K), \bX_{\partial I}\cup  (A'\times K))\stackrel{\h}{\simeq} (\bX_c\times I, \bX_c\times\partial 
I\cup A'\times I)$
 and the latter is just the product of pairs
 $(\bX_c, A')\times (I, \partial I)$. Our claim follows.

 Next, by examining the exact sequence of the triple $(X_\cD, X_S, X_c)$ and by using the annulation of $H_q(X_S, X_c)$ proved above, we see that $(X_\cD, X_c) \hookrightarrow (X_\cD, X_c)$ gives,  in homology, an isomorphism in dimensions  $q \le k+1$ and an epimorphism in $q = k+2$.
 To end our proof, we just combine this with the isomorphism
 $H_*(X_\cD, X_S) \simeq H_*(X, X_K)$, obtained by excision. 
 
\end{proof}

Since the kernel of the map $H_{k+1}(X_c) \to H_{k+1} (X)$ is equal to the image of the boundary map $H_{k+2} (X, X_c)\stackrel{\partial}{\to} H_{k+1} (X_c)$, we focus on the latter.
 Consider the commutative diagram:

\begin{equation}\label{eq:2}
\begin{array}{ccc}
 H_{k+2} (X_\cD, X_c) & \stackrel{\mbox{\tiny $\iota_*$}}{\longrightarrow} & H_{k+2} (X, X_c) \\
 \ \ \ \ \  \ \ \ \ \ \ \mbox{\tiny{$\partial_1$}} \searrow & \  & \swarrow  \mbox{\tiny{$\partial$}} \ \ \ \ \  \ \ \ \ \ \ \\
  \  & H_{k+1} (X_c) & \  \\
\end{array}
 \end{equation}
 where $\partial$ and $\partial_1$ are boundary morphisms.
  Since Proposition \ref{p:1} shows that $\iota_*$ is an epimorphism, we get:

\begin{corollary}\label{c:1} 
 If $H_q(X_c, A')=0$ for $q\le k$ then, in the diagram \rm (\ref{eq:2}), \it we have $\im \partial = \im \partial_1$.
 \fin
\end{corollary}


 Notice that, for any $M\subset \bP^1$, $X_M$ is homotopy equivalent to $\bX_M$ to which one attaches, along $A'\times M$, the product $A'\times \Cone (M)$. Since $\cD$ is contractible, it follows that $X_\cD \stackrel{\h}{\simeq} \bX_\cD$. Hence the pair $(X_\cD, X_c)$ is homotopy equivalent to $(\bX_\cD, X_c)$ and we may identify the boundary morphism
 $H_{k+2} (X_\cD, X_c) \stackrel{\mbox{\tiny{$\partial_1$}}}{\rightarrow}  H_{k+1}(X_c)$ to the boundary morphism $H_{k+2} (\bX_\cD, X_c) \stackrel{\mbox{\tiny{$\partial_1$}}}{\rightarrow}  H_{k+1}(X_c)$.
 
 Remark also that we have the excision $H_* (\cup_i \bX_{D_i}, \cup_i X_{c_i}) \stackrel{\simeq}{\rightarrow} H_* (\bX_\cD, X_c)$ which gives a decomposition of the homology $H_*(\bX_\cD, X_c)$ into the direct sum $\oplus_i H_* (\bX_{D_i}, X_{c_i})$. Then the boundary map $\partial_1$ is identified to the boundary map $\partial_2$ obtained as sum of the boundary maps $\partial_i : H_{k+2} (\bX_{D_i}, X_{c_i}) \to  H_{k+1} (X_{c_i})$, where $X_{c_i}$ is identified with $X_c$ by parallel transport along the path $\gamma_i$. 
 
 With these identifications, we have the following commutative diagram:
 
\[
\begin{array}{cccc}
 \oplus_i H_{k+2} (\bX_{D_i}, X_{c_i}) \simeq & H_{k+2} (\cup_i \bX_{D_i}, \cup_i X_{c_i}) & \  & \  \\
 \  & \mbox{\tiny{\rm exc}} \downarrow \mbox{\tiny{$\simeq$}} &  \ &  \searrow \mbox{\tiny{$\partial_2$}} \ \ \  \ \ \ \  \ \ \ \ \  \ \ \ \ \ \\
 \ & H_{k+2} (\bX_{\cD}, X_c) & \stackrel{\mbox{\tiny{$\partial_1$}}}{\longrightarrow} & H_{k+1} (X_c).
 \end{array}
 \]
 
 It then follows that
 
 \begin{equation}\label{eq:partial}
 \im \partial_1 = \sum_i \im \partial_i.
\end{equation}
 
 
 Our theorem will be proved if we do the following:\\
 
 \noindent
(i). \ Prove that $H_q(\bX_{D_i}, X_{c_i})=0$, for  $q\le k+1$ and all
$i$.

\noindent 
(ii). \  Find the image of the map $\partial_i : H_{k+2} (\bX_{D_i}, X_{c_i}) \to H_{k+1} (X_{c_i})$, for all
$i$.\\

We shall reduce these problems again, by replacing $\bX_{D_i}$ by $\bX^*_{D_i}:= \bX_{D_i}\m \Sing_\cS p$. For this, we use condition (C3) for (ii), respectively condition (C3i) for (i).

\begin{lemma}\label{l:2}
If $\Hd_{\bX \cap\Sing_\cS p} \bX \ge s+1$ then, for all $i$, the map induced by inclusion  
$H_q(\bX^*_{D_i}, X_{c_i}) \stackrel{j_*}{\rightarrow} H_q(\bX_{D_i}, X_{c_i})$ is an isomorphism, for $q \le s-1$, and an epimorphism, for $q = s$.
\end{lemma}
 \begin{proof}
 Due to the exact sequence of the triple $(\bX_{D_i},\bX^*_{D_i}, X_{c_i})$, it will be sufficient to prove, for all $i$, that $H_q(\bX_{D_i},\bX^*_{D_i})=0$, for $q\le s$. This is true since the  inclusion: 
 \[ (\bX_{D_i}\cap (\cup_j B_{ij}), \bX_{D_i}\cap  (\cup_j B_{ij}\setminus \{a_{ij}\})) \hookrightarrow (\bX_{D_i},\bX^*_{D_i})
 \]
 is an excision in homology (notice that the unions are disjoint).  As usual, $B_{ij}\subset \bX$ denotes a Milnor ball centered at the singular point $a_{ij}\in \Sing_\cS p$. 
 
 Indeed, the hypothesis $\Hd_{\bX \cap\Sing_\cS p} \bX \ge s+1$ tells that the homology of each pair $(\bX_{D_i}\cap B_{ij}, X_{D_i}\cap  B_{ij}\setminus \{a_{ij}\})$ annulates up to dimension $s$.
 \end{proof}  


\begin{corollary}\label{c:2}
If $\Hd_{\bX \cap\Sing_\cS p} \bX \ge k+3$, then, for all $i$:
\[ \im (\partial_i : H_{k+2}(\bX_{D_i}, X_{c_i}) \to H_{k+1}(X_{c_i})) = \im (\partial'_i : H_{k+2}(\bX^*_{D_i}, X_{c_i}) \to H_{k+1}(X_{c_i})).
\]
\end{corollary}
\begin{proof}  We have that $\partial'_i = \partial_i \circ j_* $, where $j_* : H_{k+2}(\bX^*_{D_i}, X_{c_i}) \to H_{k+2}(\bX_{D_i}, X_{c_i})$ is induced by the inclusion. By Lemma \ref{l:2}, $j_*$ is surjective, hence $\im \partial'_i = \im \partial_i$.
\end{proof} 

 The last step in the proof of Theorem \ref{t:main} is the following result, where the variation maps come in:

\begin{proposition}\label{p:3}
If $H_q(X_{c_i}, X^*_{a_i})=0$, for $q\le k$, then:
\begin{enumerate}
\rm \item \it $H_q(\bX^*_{D_i}, X_{c_i})=0$ for $q\le k+1$.
\rm \item \it $\im \partial'_i = \im (\var_i : H_{k+1}(X_{c_i}, X^*_{a_i}) \to H_{k+1}(X_{c_i}))$.
\end{enumerate}
\end{proposition}
\begin{proof}
Let us take Milnor data $(B_{ij}, D_i)$ at the (stratified) singularities $a_{ij}$. Recall that the radius of $D_i$ is very small in comparison to the radius of $B_{ij}$. We shall give the proof for a fixed index $i$ and therefore we suppress the lower indices $i$ in the following.

\noindent
(a). Let $D^* = D\m \{ a\}$. By retraction, we identify $D^*$ to a circle and cover this circle with the union of two arcs $I\cup J$, as follows: for the standard circle $S^1$, we take $I:= \{ \exp{i\pi t} \mid t\in [-\frac{1}{2}, 1]\}$, $J:= \{ \exp{i\pi t} \mid t\in [\frac{1}{2}, 2]\}$. Then $\bX_{D^*} \stackrel{\h}{\simeq} \bX_I \cup \bX_J$ and 
$X_c \stackrel{\h}{\simeq} \bX_J \simeq X_c\times J$. With these notations, we have the following isomorphisms induced by homotopy equivalences:
\[ H_*(\bX^*_D , X_c) \simeq H_*(\bX_{D^*} \cup X^*_a \times D, X_c\cup X^*_a \times D) \simeq H_*(\bX_I \cup \bX_J \cup X^*_a \times D, \bX_J \cup X^*_a \times D),\]
where $X^*_a \times D$ is a notation for $\bX_D \setminus \cup_j B_{ij}$, which is the total space of a trivial fibration over $D$, with fiber $X_a\setminus \cup_j B_{ij} \stackrel{\h}{\simeq} X^*_a$. 

We then excise $\bX_J \cup X^*_a \times D$ from the last pair and get the homology of the pair $(\bX_I, X_{\partial I} \cup X^*_a \times I)$, which pair is homotopy equivalent to the product $(X_c, X^*_a) \times (I, \partial I)$.
Since, by hypothesis, the homology of the pair $(X_c, X^*_a)$ annulates up to dimension $k$, it follows that the homology of the last product annulates up to dimension $k+1$.

\noindent
(b). In the following commutative diagram, the variation map identifies to the right-hand vertical arrow. This diagram is a Wang type exact sequence, the proof of which is explained by Milnor in \cite[page 67, Lemma 8.4]{Mi}.

\[ \begin{array}{ccc}
H_{k+2} (\bX^*_{D_i} , X_{c_i}) & \stackrel{\partial'_i}{\longrightarrow} & H_{k+1} (X_{c_i}) \\
\mbox{\tiny excision} \uparrow \mbox{\tiny $\simeq$} & \  & \uparrow \mbox{\tiny $\var_i$} \\
 H_{k+2} (\bX_I, X_{\partial I} \cup X^*_{a_i} \times I) & \simeq & H_{k+1} (X_{c_i}, X^*_{a_i})\otimes H_1(I, \partial I)
\end{array} \]

This shows that $\im \partial'_i = \im \var_i$.

\end{proof}

We are now able to conclude the proof of Theorem \ref{t:main}.
The claim (i) above, and hence the first claim of the theorem, follows from Lemma \ref{l:2} and Proposition \ref{p:3}(a).

The second claim of the theorem follows by the sequence of results:
Corollary  \ref{c:1}, equality (\ref{eq:partial}), Corollary  \ref{c:2} and Proposition \ref{p:3}(b).

\section{Further results and particular cases}\label{s:comments}

\subsection{The case $A\subset V$}\label{ss:subset}\ 
We discuss in the following the case $A'=\emptyset$, equivalently,
$A\subset V$, which is complementary to 
the one we have considered until now.  One would be tempted to replace  the condition (C1)  with ``$H_q(X_c) =0$, for $q\le k$", but this appears to be too restricting.

Nevertheless, in case $h_{|X}$ is not onto $\bP^1$, the situation becomes
more interesting. So let us assume that $V$ contains a
fiber of the pencil $h : Y\m A \to \bP^1$. Even if the axis $A$ is outside the
space $X$, the ``singularities in the axis" influence the topology of the pencil.
We have the following result on a class of nongeneric pencils,
 disjoint from the class considered in Theorem \ref{t:main}.
 Let us denote $\Sigma := \sigma (\Sing_\cS p)$.

\begin{theorem}\label{t:2} 
  Let $X= Y\m V$ have a structure of Lefschetz fibration with isolated
  singularities, such that $V$ contains a member of the pencil. For some fixed $k \ge 0$, assume that $H_q(X_c, X^*_{a_i})=0$ for $q\le k$, where $X_c$ is a general member $X_c$ and $X_{a_i}$ is any atypical one. We have:
  
\begin{enumerate}
\rm \item \it If $H_q(X, X\m \Sigma)=0$ for $q\le k+1$, then 
$H_q(X,X_c)=0$ for $q\le k+1$.
\rm \item \it If $H_q(X, X\m \Sigma)=0$ for $q\le k+2$, then: 
\[  H_{k+1}(X) \simeq  H_{k+1}(X_c)/\sum_i^p \im \var_i. \] 
\end{enumerate} 
\end{theorem}

\begin{proof}
 The proof follows the lines of the proof of Theorem \ref{t:main} and we shall only point out the differences, using the same notations. In our case, the target of the holomorphic function $h_{|X}$ is $\bP^1\m \{\alpha\}$ for some $\alpha\in \bP^1$.  We have $\cD \stackrel{\h}{\simeq}\bP^1\m \{\alpha\}$ and therefore $X_\cD\stackrel{\h}{\simeq} X$. Examining the proofs of Proposition \ref{p:1} and Corollary \ref{c:1}, we see that, under our assumptions, their conclusions hold without any restrictions on $k$. Hence (C1)  does not enter as condition in our proof. On the other hand, by Proposition \ref{o:1}(b), we can use (C3)' instead of (C3). Condition (C2) is itself an assumption of the above theorem.
\end{proof} 

\subsection{Comparing to the $\rHd$ condition}\label{ss:rHd}\

\begin{proposition}\label{n:cond}
Theorem \ref{t:main} holds if we replace the conditions {\rm (C2)} and {\rm (C3)} by the single condition:

\smallskip\smallskip
\noindent
{\bf (C4)} \ \ \ $\rHd X \ge k+3$.

\smallskip\smallskip

The first claim of Theorem \ref{t:main} holds with a weaker assumption in place of {\rm (C4)}, namely (see Note \ref{n:main}):

\smallskip\smallskip
\noindent
{\bf (C4i)}\ \ \ $\rHd X \ge k+2$.


\end{proposition}
\begin{proof} Indeed, $\rHd X \ge q$ implies $\rHd \bX \ge q$, since $\bX$ is a hypersurface in $X\times \bP^1$ and one can apply the result of Hamm and L\^e \cite[Theorem 3.2.1]{HL}. This in turn implies 
$\Hd_\alpha \bX \ge q$, for any point $\alpha\in \bX$, by definition.

Next, $\rHd \bX \ge q$ implies that the homology of the pair $(\bX_{D_i}, X_{c_i})$ annulates up to dimension $q-1$, by \cite[Proposition 3.4]{Ti-lef} (where $\rhd$ is used instead of $\rHd$, but the proof is the same). This shows that conditions 
(C1) + (C4i) imply the first claim of Theorem \ref{t:main}.

Furthermore, if we assume (C4) instead of (C4i), then, besides the annulation of the homology of $(\bX_{D_i}, X_{c_i})$ up to $k+2$ (shown just above), it follows that $H_q(\bX^*_{D_i}, X_{c_i})=0$ for $q\le k+1$,  by Lemma \ref{l:2}.
The proof of Proposition \ref{p:3}(a) shows in fact that the annulation of homology of $(\bX^*_{D_i}, X_{c_i})$ up to $k+1$ is equivalent to the annulation of homology of the pair $(X_{c_i}, X^*_{a_i})$ up to $k$, which is condition (C2). Now Theorem \ref{t:main} applies.
\end{proof}

\subsection{Particular cases}\label{ss:relax}\
From Theorem \ref{t:main} and its proof, one may derive several versions in particular cases, recovering some of the results in the literature. To do that, one has to take into account the following observations (still under the condition $A\cap X \not= \emptyset$):
 
\begin{proposition}\label{o:1}\hspace*{\fill}
\begin{enumerate}
\item In case $\bX\cap \Sing_\cS p = \emptyset$, the condition {\rm (C3)} is void.
\item In case $(A\times \bP^1)\cap\bX \cap \Sing_\cS p =\emptyset$, we may replace condition {\rm (C3)} by the following more general condition (which is also more global):

\vspace{2mm}
\noindent
\hspace*{-8mm} {\bf (C3)'} \ \    $H_q(X, X\m \Sigma) =0$, for  $q \le k+2$.
\vspace{2mm}
\item In case $(A\times \bP^1)\cap\Sing_\cS p =\emptyset$, if condition {\rm (C1)} is true, then {\rm (C2)}
is equivalent to the following:

\vspace{2mm}
\noindent
\hspace*{-8mm} {\bf(C2)'} \ \ \   $H_q (X^*_{a_i}, X^*_{a_i}\cap A) =0$,  for  $q \le k-1$.

\end{enumerate}
\end{proposition}
\begin{proof}
(a). is obvious. \\
(b). By examining the Proof of Theorem \ref{t:main}, we see that we have used the homology depth condition only to compare $\bX_{D_i}$ to $\bX^*_{D_i}$. We may cut off from the proof this comparison (which means Lemma \ref{l:2} and Corollary \ref{c:2}) and start from the beginning with the space $X\m \Sigma$ instead of the space $X$. Taking into account that, under our hypothesis, $\bX^*_{D_i} = \bX^*_{D_i}\m \Sigma$, for all $i$, the effect of this change is that
the proof yields the conclusion ``$H_q(X\m\Sigma , X_c)=0$, for $q\le k+1$'' and the corresponding statement for the vanishing cycles. At this final stage, condition (C3)' allows one to replace $X\m\Sigma$ by $X$. 

\noindent
(c) When there are no singularities in the axis, we have $A\cap X^*_{a_i} = A\cap X_c$, for any $i$. Then the exact sequence of the triple 
$(X_c, X^*_{a_i}, A\cap X^*_{a_i})$ shows that the boundary morphism
\[  H_q(X_c, X^*_{a_i}) \to H_{q-1}(X_c, A\cap X^*_{a_i})\]
is an isomorphism, for $q\le k$, by condition (C1). This implies our claimed equivalence.
\end{proof}
 
In case of quasi-projective varieties, we have an abundance of hyperplane pencils, which are moreover generic, in the sense that
the axis is transversal to the stratification. It easily follows that such a pencil has no singularities along the axis (see e.g. the proof of \cite[Proposition 2.4]{Ti-lef} for a detailed explanation). We are therefore in the conditions of Proposition \ref{o:1}(b) and (c).
Another nice aspect of quasi-projective varieties is that the Lefschetz structure is hereditary on slices. Namely, as already observed by Ch\'eniot \cite{Ch1}, since the axis $A$ is chosen to be generic, it becomes in turn a generic slice of a hyperplane slice of $X$, and so on.

Condition (C2)' has been used by Ch\'eniot \cite{Ch1, Ch} and Eyral \cite{Ey} in theorems on generic pencils of hyperplanes,  respectively condition (C3)' has been used by C. Eyral in proving a version of the First Lefschetz Hyperplane Theorem (compare to \cite[Proof of Theorem 2.5]{Ey}). Therefore, via Proposition \ref{o:1} and the preceeding observations in case of quasi-projective varieties, our Theorem \ref{t:main}
 recovers the results in the cited articles. 
%
\section{Vanishing cycles of polynomial functions on $\bC^n$} \label{s:poly}

A polynomial function $f: \bC^n \to \bC$ can naturally be considered as a nongeneric pencil of hypersurfaces on $\bC^n$, which is a particular quasi-projective variety.
Indeed, this function extends as a meromorphic function on $\bP^n$, as follows.
 If $\deg f =d$, then $h = \tilde f/z^d : \bP^n \dashrightarrow \bP^1$, where $\tilde f$ is the homogenized of $f$ with respect to the new variable $z$ and the axis of the pencil is $A = \{ f_d =0\} \subset H_\ity$. Here we have $Y=\bP^n$, $V=H_\ity = \{ z=0\} \subset \bP^n$. We are in the situation described in \S \ref{ss:subset}, namely we have a pencil  on $X= \bC^n$, where $h_{|\bC^n} = f$. In particular, $\Sigma = \Sing f$. 
 
 For such a pencil, we may work under more general hypotheses: we assume that the function $f$ has isolated singularities, but we put {\em no condition on singularities in the axis}, which may be non-isolated. 
 We show how this can fit in the theory developed before.
 
 Take 
 the complement of a big ball $B\subset \bC^n$, centered at the origin of a fixed system of coordinates on $\bC^n$. The complement $C_B :=\bC^n \m B$ plays the role of a ``uniform" neighbourhood of the whole hyperplane at infinity $H_\ity$ and of all singularities in the axis together. For big enough radius of $B$, we have  
 \[X_{a_i} \cap B \stackrel{\h}{\simeq} X_{a_i},\]
  for any $i$, since the distance function has a finite set of critical values on the algebraic sets $X_{a_i}$. We claim that $f^{-1}(D_i) \cap B \setminus \cup_j B_{ij} \to D_i$ is a trivial fibration, where the $B_{ij}$'s are small Milnor balls around the critical points of $f$ on $X_{a_i}$ and $D_i$ is a small enough disk. Indeed, the fibers of $f$ over $D_i$ are transversal to the boundary of a big ball and transversal to the boundaries of the Milnor balls. Our claim then follows by Ehresmann's Theorem.
  
    This implies, as in \S \ref{piece:var}, that there is a well defined geometric monodromy representation at each $a_i\in \Lambda\subset \bC$, $\rho_i : \pi_i(\partial\bar D_i , c_i) \to \Iso(X_{c_i}, X_{c_i}\m (C_B\cup\cup_j B_{ij}))$. This induces a variation map:
  \[ \var_i : H_k (X_{c_i}, X^*_{a_i}) \to  H_k (X_{c_i}),\]
 where $X^*_{a_i} := X_{a_i} \m \Sing f$ is now used as a notation for the subset $X_{c_i}\m (C_B\cup\cup_j B_{ij})$ of $X_{c_i}$. This is justified by the fact that $X_{a_i} \m \Sing f$ is homotopy equivalent to 
 $X_{a_i}\cap B\setminus \cup_j B_{ij}$, which in turn can be identified to $X_{c_i}\m (C_B\cup\cup_j B_{ij})$ as fibers in the above mentioned trivial fibration.
 
 We shall show that Theorem \ref{t:2} holds for a pencil
 defined by a polynomial function with isolated singularities $f :\bC^n \to \bC$, without any condition on singularities at the axis at infinity and   
 moreover, that we have a more precise grip on variation maps. 
 
 Let us first remark that the boundary map $H_{*+1}(\bC^n, X_c) \stackrel{\partial}{\to} \tilde H_{*}(X_c)$ is an isomorphism in any dimension. This follows from the long exact sequence of the pair $(\bC^n, X_c)$.
 
Next, we have by excision: $H_{*+1}(\bC^n, X_c) \simeq \oplus_i H_{*+1}(X_{D_i}, X_c)$. These show that $H_{*}(X_c)$ decomposes into the direct sum of vanishing cycles at each atypical fiber $X_{a_i}$.
Note that the direct sum decomposition depends on the paths $\gamma_i$.
  
  We say that $\im ( H_{*+1}(X_{D_i}, X_c) \stackrel{\partial_i}{\to} H_*(X_c))$ is the module of {\em vanishing cycles} at the fiber $X_{a_i}$.
  It has been proved in general that $H_{*}(X_c)$ is the direct sum of the modules of vanishing cycles (see \cite[proof of Theorem 3.1]{ST}, \cite[Theorem 1.4]{NN1}) regardless of the singularities of $f$. 
 
  It is well-known that in case of a holomorphic function germ with isolated singularity
 on $\bC^n$, the variation map of the local monodromy
 is an isomorphism \cite{Mi}. But in our global case of a polynomial function  with isolated singularities, the variation maps cannot be isomorphisms since the homology of the fiber $H_*(X_c)$ captures information on vanishing cycles at all atypical fibers $X_{a_i}$ together. We may prove
 the following statement, the part (b) of which being just Neumann-Norbury's
 result \cite[Theorem 2.3]{NN2} via an identification (by some excision) of our variation map to the local variation maps used in \cite{NN2}.      

\begin{proposition} \label{t:poly}
Let $f :\bC^n \to \bC$ be a polynomial function with isolated singularities.
Then:
\begin{enumerate}
\rm \item \it If $H_q(X_c, X^*_{a_i})=0$ for $q\le k$ and for any $i$, then $\tilde H_q(X_c)=0$ for $q\le k$.
\rm \item \cite[Theorem 2.3]{NN2}\ \it 
The variation map $\var_i : H_*(X_{c_i}, X^*_{a_i}) \to  
H_*(X_{c_i})$ is injective, for any $i$. In particular, we have $H_q(X_{c}) \simeq \sum_i \im(\var_i)$ for the first integer $q\ge 1$ such that $H_q(X_c) \not= 0$.
\end{enumerate}
\end{proposition}
\begin{proof}
Since the fibers of $f$ are Stein spaces of dimension $n-1$, their homology groups are trivial in dimensions $\ge n$. The condition (C3)' is largely satisfied, since $(\bC^n, \bC^n \m \Sing f)$ is $(2n-1)$-connected.
Hence part (a) follows from Theorem \ref{t:2}.
For part (b), remark first that, by the above arguments, the boundary map
$\partial_i : H_{*+1}(X_{D_i}, X_{c_i}){\to} \tilde H_{*}(X_{c_i})$ is injective, for any $i$. Next, one may replace $X_{D_i}$ by $X^*_{D_i}$ since
$(X_{D_i}, X^*_{D_i})$ is $(2n-1)$-connected. It follows that the boundary
morphism $\partial'_i : H_{*+1}(X^*_{D_i}, X_{c_i}){\to} \tilde H_{*}(X_{c_i})$ is injective. As in Proposition \ref{p:3}, one may identify $H_{*+1}(X^*_{D_i}, X_{c_i})$ to $H_*(X_{c_i}, X^*_{a_i})$, by excision
 and $\partial'_i$ can be identified with 
$\var_i : H_*(X_{c_i}, X^*_{a_i}) \to  
H_*(X_{c_i})$.
\end{proof}

The image of the ``pseudo-embedding'' $\iota : X^*_{a_i} \stackrel{\h}{\simeq} X_{a_i}\cap B\setminus \cup_j B_{ij} \hookrightarrow X_{c_i}$ plays here the role of the boundary of the Milnor fiber in the local case. We may therefore call $\im \iota_*$ the group of ``boundary cycles" at $a_i$. We immediately get the following consequence; it can also be deduced, by a series of identifications, from the Neumann-Norbury more general result  \cite[Theorem 1.4]{NN1}.

\begin{corollary}\label{c:poly}
The invariant cycles under the monodromy at $a_i$ are exactly the boundary cycles, i.e. the following sequence is exact:
 \[ H_*(X^*_{a_i}) \stackrel{\iota_*}{\to} H_*(X_{c_i}) \stackrel{\nu_a - \id}{\to} H_*(X_{c_i}).\]
\end{corollary}
\begin{proof} 
We have the following commutative diagram, where the first row is 
the exact sequence of the pair $(X_{c_i}, X^*_{a_i})$:
\[ \begin{array}{ccccc}
 H_*(X^*_{a_i}) & \stackrel{\iota_*}{\longrightarrow} & H_*(X_{c_i}) & \stackrel{j_*}{\longrightarrow} & H_*(X_{c_i}, X^*_{a_i}) \\
 \ & \ & \mbox{\tiny $\nu_i - \id$} \searrow &  \ & \swarrow \mbox{\tiny $\var_i$}  \\
 \ & \ & \ & \  H_*(X_{c_i})
 \end{array}
 \]
 We have that $\im \iota_* = \ker j_*$.
 Since $\nu_i -\id = \var_i \circ j_*$, and since $\var_i$ is injective
 by Proposition \ref{t:poly}, our claim follows.
\end{proof}

\begin{note}\label{n:clemens}
This result may be considered as a counterpart, in a non-proper situation, of the well-known ``invariant cycle theorem'' proved by Clemens \cite{Cl}. The later holds for proper holomorphic functions $g: X\to D$, in cohomology (thus ``invariant co-cycle theorem'' would be more appropriate), where $X$ is a K\"ahler manifold. It sais that the following sequence is exact:
$ H^*(X) \stackrel{\i^*}{\to} H^*(X_c)  \stackrel{h - \id}{\to} H^*(X_c),$
 where $h$ denotes the monodromy around the center of the disk $D$ (assumed to be the single critical value of $g$).
 
  Is is natural to ask if an invariant cycle result similar to Corollary \ref{c:poly} holds for more general classes of non-proper pencils.
\end{note}


\end{document}